\documentclass[11pt,draft]{article}
\usepackage{amsfonts}
\usepackage{mathrsfs}
\usepackage{amsmath,amsfonts,mathrsfs,amssymb,color}
\usepackage{indentfirst}

\numberwithin{equation}{section}

\newtheorem{theo.}{\quad\, Theorem}[section]
\newtheorem{defi.}{\quad\, Definition}[section]
\newtheorem{lemm.}{\quad\, Lemma}[section]
\newtheorem{coro.}{\quad\, Corollary}[section]

\begin{document}

\title{Positive solutions for nonlocal extended Fisher-Kolmogorov and
Swift-Hohenberg equations via bifurcation methods$^*$ }
\author{Jinxiang Wang$^{\ast}$
\\
 {\small  Department of Applied Mathematics, Lanzhou University of Technology, Lanzhou, P.R. Chin}\\
}
\date{} \maketitle
\footnote[0]{E-mail address: wjx19860420@163.com(Jinxiang Wang), \ \ Tel:
86-931-7971297.
} \footnote[0] {$^*$Corresponding author: Jinxiang Wang.\  Supported by the NSFC(No.11801453). }

 \begin{abstract}
\baselineskip 18pt
In this paper we study the existence of positive solutions for a class of nonlocal fourth-order nonautonomous differential equations which can be seen as generalization of extended Fisher-Kolmogorov and Swift-Hohenberg equations. The main result is proved using bifurcation theory.

 \end{abstract}

\vskip 3mm
{\small\bf Keywords.} {\small Fourth order boundary value
problem, Nonlocal problem, Global bifurcation, Positive solution}

\vskip 3mm

{\small\bf MR(2000)\ \ \ 34B10, \ 34B18}

\baselineskip 22pt

\section{Introduction}
\vskip 3mm
Consider the fourth-order differential equation
$$
u^{(4)}(x)-pu''(x)-au(x)+bu^{3}(x)=0,\ \ \ \ \ x\in[0,L], \eqno (1.1)
$$
subject to the boundary conditions
$$u(0)=u(L)=u''(0)=u''(L)=0,\eqno (1.2)$$
where $u=u(x)$ is a function of the space variable $x$ in $[0,L]$, $p,a,b$ are constants. In studies of pattern formation, equation (1.1) plays an important role. When $p<0$, (1.1) is called the Swift-Hohenberg equation [1,2], and for $p>0$ it is called the extended Fisher-Kolmogorov equation [3]. Problem (1.1),(1.2) occurs in a variety of applications: the behavior close to a so-called Lifshitz point in phase transition physics (e.g. for nematic liquid crystals and for erroelectric crystals) [4], the rolls in a Rayleigh-Benard convection cell (two parallel plates at different temperatures with a liquid in between) [5], the waves on a suspension bridge [6,7], the buckling of a strut on a nonlinear elastic foundation [8], and pulse propagation in optical fibers [9] and for other references see [10].

Solutions of (1.1) which are bounded on the real line have been studied by several researchers, see e.g. [11-26]. In [13], when the coefficients $a$ and $b$ are even and $2L$ positive continuous periodic functions, Chaparova established a multiplicity result of $2L$ periodic solutions for the problem (1.1),(1.2) using variational techniques. In [27], taking into account  the interaction-induced modification of the environment around the
individual, Ma and Dai discussed the spatially nonlocal generalization of (1.1), that is the equation
$$
u^{(4)}(x)-pu''(x)-a(x)u(x)+u(x)\int_0^Lf(|x-y|)u^{2}(y)dy=0,\ \ \ \ x\in(0,L),\eqno (1.3)
$$
where the nontrivial integral kernel $f(|x-y|): \ [0,L]\times [0,L] \rightarrow \mathbb{R}^{+}$ is a continuous and monotone decreasing function.  By using Clark¡¯s theorem and symmetric mountain-pass
theorem, existence of nontrivial periodic solutions for problem (1.3),(1.2) is proved in [27].  For further references about the solution of fourth order problems like (1.1),(1.2) which have broad classes of nonlinearities $g(u)$ or $g$  depending on $u$ and its derivatives,  we refer to [28-43].

 Motivated by the above works described, in this paper we are going to study the following nonlocal problem
 $$
   \left\{\begin{array}{ll}
     u^{(4)}(x)-p(x)u''(x)-a(x)u(x)+u^{\rho}(x)\int_0^1f(x,y)u^{\sigma}(y)dy=0,\ \ \ \ \ \ \ x\in(0,1),\\
    u(0)=u(1)=u''(0)=u''(1)=0\\
\end{array} \right.\eqno (1.4)
$$
where $p$ and $a$ are continuous functions, $\rho\geq 1$ and $\sigma>0$ are constants, $f\in L^{\infty}([0,1]\times[0,1])$ is a nonnegative function satisfying other hypotheses that will be detailed below. When the variable coefficient $p$ is identically equal to a constant, $\rho\equiv 1, \sigma\equiv2$ and $f$ is as in (1.3),  it is easy to see that problem (1.4) will degenerate into problem (1.3),(1.2). Unlike the works in [27], in this note we are going to apply global bifurcation theory in order to study the existence of positive solutions for (1.4). 
It should be pointed that,  global bifurcation phenomena for fourth order problems like (1.1),(1.2) with different classes of nonlinearities have been investigated in [40-43], 
but as far as we know, there have been no studies on the bifurcation phenomena for nonlocal fourth order problems.

We first introduce the class $K$, which is formed by functions $f: [0,1]\times[0,1] \rightarrow \mathbb{R}$ verifying:

(i) $f\in L^{\infty}([0,1]\times[0,1])$ and $f(x,y)\geq0$ for all $x,y\in [0,1]$.

(ii) If $w\in C[0,1]$ and $\int_0^1\int_0^1f(x,y)|w(y)|^{\sigma}w^{\rho+1}(x)dxdy=0$, then $w\equiv0$.

Our main results are:

\noindent{\bf Theorem 1.1 } Suppose that $f\in K$. Assume the variable coefficients $a,b$ satisfying\\
\noindent(H1) $p\in C[0,1]$ with $p(x)>-\pi^{2}, \ x\in [0,1]$;\\
\noindent(H2) $a\in C[0,1]$ with $a\geq0$ on $[0,1]$ and $a \not\equiv 0$ on any subinterval of $[0,1]$;\\
\noindent(H3) $$\pi^{4}+2\pi^{2}\int_{0}^{1}p\sin^{2} (\pi x)dx<2\int_{0}^{1}a\sin^{2} (\pi x)dx, $$
then problem (1.4) has a positive solution.

The rest paper is arranged as follows:  In Section 2, we give some preliminaries and show a global bifurcation
phenomena of the corresponding auxiliary problem with parameter. In Section 3, we discuss the direction of the component and complete the proof of Theorem 1.1.

\vskip 3mm

\section{Preliminaries and the Global bifurcation}

Let $X=C[0,1]$ be the Banach space of continuous function defined on $[0,1]$, with its usual normal $||\cdot||_{\infty}$, it is easy to see that $P:=\{u\in C[0,1]: u(x)\geq 0, \forall\  x\in[0,1]\}$ is a positive cone in $X$. Let $E:=\{u\in C^{2}[0,1]:u(0)=u(1)=u''(0)=u''(1)=0\}$ with the norm $\|u\|_{E}=\max\{\|u\|_{\infty},\|u'\|_{\infty},\|u''\|_{\infty}\}$. 

Define an linear operator $L: C^{4}[0,1]\cap E \to X$
$$
Lu:=u''''-p(x)u''.
$$  Applying  Elias's theory, Ma [40] has verified that if (H1) hold, 
then $L$ is a disconjugate, and consequently, $L$ is positive, and $L^{-1}: E\rightarrow E $ is completely continuous. Moreover, by using Elias¡¯s eigenvalue theory [44], Ma [40] also proved that\\
\noindent{\bf Lemma 2.1} \ Assume (H1) and (H2) hold, then the eigenvalue problem
$$
   \left\{\begin{array}{ll}
     u^{(4)}(x)-p(x)u''(x)=\lambda a(x)u(x),\ \ \ \ \ \ \ x\in(0,1),\\
    u(0)=u(1)=u''(0)=u''(1)=0.\\
\end{array} \right.\eqno (2.1)
$$ has an infinite sequence of positive eigenvalues $$\lambda_{1}<\cdots<\lambda_{k}<\cdots,\eqno (2.2)$$ and to each eigenvalue $\lambda_{k}$ there exists an essential unique eigenfunction $\phi_{k}$ which
has exactly $k-1$ simple zeros in (0,1). In particular, the eigenfunction $\phi_{1}$ corresponding to the principle eigenvalue $\lambda_{1}$ is positive.

It is easy to see that there exist $C>0$ such that for each $g\in C[0,1]$, there exist a unique $v\in E$ satisfying
 $$
 Lv=g(x),\ \ \ \ \ \ \ x\in(0,1),\eqno (2.3)
$$
and
$$\|v\|_{E}\leq C\|g\|_{\infty}.\eqno (2.4) $$

Assuming that $f\in K$, for any $w\in L^{\infty}(0,1)$, we will consider the function $\theta_{w}: [0,1]\rightarrow \mathbb{R}$ given by
$$ \theta_{w}(x)=\int_0^1f(x,y)|w(y)|^{\sigma}dy.$$ Once that $f$ and $w$ are bounded, we have that $\theta_{w}$ is well defined. 

Using the above notation, it is easy to observe that $u$ is a solution of (1.4) if and only if $u$ is a positive solution of
$$
   \left\{\begin{array}{ll}
     u^{(4)}(x)-p(x)u''(x)=a(x)u(x)-u^{\rho}(x)\theta_{u},\ \ \ \ \ \ \ x\in(0,1),\\
    u(0)=u(1)=u''(0)=u''(1)=0\\
\end{array} \right.
$$

Let us consider
$$
Lu=\lambda a(x)u(x)-\lambda u^{\rho}(x)\theta_{u},\eqno (2.5)
$$
where $\lambda\geq0$ is parameter. Obviously, (2.5) has trivial solution $u\equiv 0$ for any $\lambda$. In what follows, we will investigate the bifurcation from the trivial solution.
Equation (2.5) can be converted to the equivalent equation
$$
u=\lambda L^{-1} (a(x)u(x))-\lambda L^{-1} (u^{\rho}(x)\theta_{u}). \eqno (2.6)
$$
Since for all $u\in C[0,1] $, $$
||\theta_{u}||_{\infty}\leq ||f||_{\infty}||u||_{\infty}^{\sigma},\eqno (2.7)
$$
this combining with (2.4) and (2.7) conclude that  $$
\|L^{-1} (u^{\rho}(x)\theta_{u})\|_{E}\leq C ||\theta_{u}||_{\infty}\|u\|_{\infty}^{\rho}\leq C ||f||_{\infty}||u||_{\infty}^{\sigma+\rho}.\eqno (2.8)
$$
we have
$$
\frac{\|L^{-1} (u^{\rho}(x)\theta_{u})\|_{E}}{||u||_{E}}\leq\frac{C ||f||_{\infty}||u||_{\infty}^{\sigma+\rho}}{||u||_{\infty}},
$$
Since $\sigma+\rho>1 (\rho>1, \sigma>0)$, then $$
\lim\limits_{u\rightarrow 0}\frac{\|L^{-1} (u^{\rho}(x)\theta_{u})\|_{E}}{||u||_{E}}=0,
$$
that is
$$
\|L^{-1} (u^{\rho}(x)\theta_{u})\|_{E}=\circ (||u||_{E})
$$
Now, we are in the position in applying Rabinowitz's global bifurcation theorem [45] to problem (2.5), we have

\noindent{\bf Lemma 2.2.} Assume that (H1) and (H2) hold, then from $(\lambda_{1},0)$ there emanate a subcontinuum $\mathcal{C}$ of positive solutions of (2.5) in the set $\mathbb{R}\times E$, which satisfies

 (i) $\mathcal{C}$ is unbounded in the set $\mathbb{R}\times E$, or

 (ii) there exists some eigenvalue $\lambda_{k}$ of (2.1) which satisfy $\lambda_{k}\neq \lambda_{1}$ and  $(\lambda_{k},0)\in \mathcal{C}$.

Since (2.5) has only trivial solution when $\lambda=0$, then by Lemma 2.1, $\mathcal{C}$ can not passes through the hyperplane $0\times E$, that is, $\mathcal{C}\subset\mathbb{R}^{+}\times E.$

\noindent{\bf Lemma 2.3.} $\mathcal{C}$ is unbounded.

\noindent{\bf Proof. } Suppose that $\mathcal{C}$ is bounded, then the (ii) of Lemma 2.1 holds, that is $\mathcal{C}$ must contain $(\lambda_{k},0)$, where $\lambda_{k}\neq \lambda_{1}$ is a eigenvalue of (2.1). Then there exist sequence $\{(\overline{\lambda}_{n},u_{n})\}\subset \mathcal{C}$ satisfying $\{u_{n}\}\in E\setminus 0$ with $\|u_{n}\|_{E}\longrightarrow 0$ and $\{\overline{\lambda}_{n}\}\rightarrow \lambda_{k}$.  Moreover, for any $n\in\mathbb{N}$, $(\overline{\lambda}_{n},u_{n})$ satisfies
$$Lu_{n}=\overline{\lambda}_{n} a(x)u_{n}(x)-\overline{\lambda}_{n} u_{n}^{\rho}(x)\theta_{u_{n}},\eqno (2.9)$$
divide (2.9) by $\|u_{n}\|_{\infty}$ and set $v_{n}=\frac{u_{n}}{\|u_{n}\|_{\infty}}$, then we get
$$Lv_{n}=\overline{\lambda}_{n} a(x)v_{n}(x)-\overline{\lambda}_{n} v_{n}(x) u_{n}^{\rho-1}(x)\theta_{u_{n}},\eqno (2.10)
$$
Then by (2.4) and (2.7), (2.10) imply that
$$
\|v_{n}\|_{E}\leq \overline{\lambda}_{n}C ||v_{n}||_{\infty}( ||a||_{\infty}+||u_{n}^{\rho-1}||_{\infty}||f||_{\infty}||u_{n}||_{\infty}^{\sigma}).
\eqno (2.11)$$
Since $\rho>1, ||v_{n}||_{\infty}=1$ and $\{u_{n}\}$ is bounded, then (2.11) imply that $\{v_{n}\}$ is bounded in $E$. Then by the Ascoli-Arzela theorem, a subsequence of $\{v_{n}\}$ uniformly converges to a limit $v\in C^{1}[0,1]$, and we again denote by $\{v_{n}\}$ the subsequence. Passing to the limit in
$$v_{n}=L^{-1}[\overline{\lambda}_{n} a(x)v_{n}(x)-\overline{\lambda}_{n} v_{n}(x) u_{n}^{\rho-1}(x)\theta_{u_{n}}],\eqno (2.12)
$$
we get
$$v=L^{-1}[\lambda_{k} a(x)v(x)],\eqno (2.13)
$$
that is, $v$ is a eigenfunction of (2.1) corresponding to the eigenvalue $\lambda_{k}$. Since $\lambda_{k}\neq \lambda_{1}$, then $v$ must change sign. Then, for $n$ large enough, $v_{n}$ must change sign, this is contradiction with $\{(u_{n},\overline{\lambda}_{n})\}\subset \mathcal{C}$.
\hfill{$\Box$}

\vskip 3mm

\section{The proof of Theorem 1.1}

\noindent{\bf Lemma 3.1.} Assume that (H1),(H2) hold. Then, $\sup\{\lambda|\ (\lambda,u)\in\mathcal{C}\}=+\infty.$

\noindent{\bf Proof.}\ \ Assume on the contrary that $\sup\{\lambda|\ (\lambda,u)\in\mathcal{C}\}=:c_{0}<\infty.$ Let $\{(\overline{\lambda}_{n},u_{n})\}\subset\mathcal{C}$ be such that $\overline{\lambda}_{n}+\|u_{n}\|_{E}\rightarrow +\infty$, then $\|u_{n}\|_{E}\rightarrow+\infty.$
Since $(\overline{\lambda}_{n},u_{n})\in\mathcal{C}$, then (2.9) hold,  multiplying (2.9) by $u_{n}$ and integrating it over $[0,1]$, based on boundary conditions and integration by parts we obtain
$$
\int_{0}^{1} [(u_{n}''(x))^{2}- p(x)u_{n}''(x)u_{n}(x)]dx=\overline{\lambda}_{n}\int_{0}^{1} [ a(x)u_{n}^{2}(x)-u_{n}^{\rho+1}(x)\int_{0}^{1}f(x,y)u_{n}^{\sigma}(y)dy]dx.\eqno (3.1)
$$
Divide (3.2) by $\|u_{n}\|_{E}^{\rho+\sigma+1}=[\max\{\|u\|_{\infty},\|u'\|_{\infty},\|u''\|_{\infty}\}]^{\rho+\sigma+1}$ , we get
$$
\frac{\int_{0}^{1} [(u_{n}''(x))^{2}- p(x)u_{n}''(x)u_{n}(x)]dx}{\|u_{n}\|_{E}^{\rho+\sigma+1}}=\frac{\overline{\lambda}_{n}\int_{0}^{1} a(x)u_{n}^{2}(x)dx}{\|u_{n}\|_{E}^{\rho+\sigma+1}}
-\frac{\overline{\lambda}_{n}\int_{0}^{1}(u_{n}^{\rho+1}(x)\int_{0}^{1}f(x,y)u_{n}^{\sigma}(y)dy)dx}{\|u_{n}\|_{E}^{\rho+\sigma+1}}
$$
Since $\overline{\lambda}_{n}\leq c_{0}$, $\rho+\sigma+1>2$ and $\|u_{n}\|_{E}\rightarrow+\infty,$ then passing to the limit in the above equality and denoting $w_{n}=\frac{u_{n}}{\|u_{n}\|_{E}}$, we can conclude that
$$
\frac{\overline{\lambda}_{n}\int_{0}^{1}(u_{n}^{\rho+1}(x)\int_{0}^{1}f(x,y)u_{n}^{\sigma}(y)dy)dx}{\|u_{n}\|_{E}^{\rho+\sigma+1}}
=\overline{\lambda}_{n}\int_{0}^{1} [w_{n}^{\rho+1}(x)\int_{0}^{1}f(x,y)w_{n}^{\sigma}(y)dy]dx\rightarrow 0.\eqno (3.2)
$$
Noticing that  $\|w_{n}\|_{E}\equiv 1$, then by Ascoli-Arzela theorem, a subsequence of $\{w_{n}\}$ uniformly converges to a limit $w\in C[0,1]$, and consequently, from Fatou Lemma, (3.2) implies that
$$
\int_{0}^{1} [w^{\rho+1}(x)\int_{0}^{1}f(x,y)w^{\sigma}(y)dy]dx\leq\lim\limits_{n\rightarrow \infty}\int_{0}^{1} [w_{n}^{\rho+1}(x)\int_{0}^{1}f(x,y)w_{n}^{\sigma}(y)dy]dx=0,\eqno (3.3)
$$
that is
$$
\int_{0}^{1} [w^{\rho+1}(x)\int_{0}^{1}f(x,y)w^{\sigma}(y)dy]dx=0,\eqno (3.4)
$$
Since $f\in K$, then (3.4) implies $w\equiv0$, that is $\{w_{n}\}$ uniformly converges to $0$ in $C[0,1]$.

On the other hand, divide (2.9) by $\|u_{n}\|_{E}$, we have
$$Lw_{n}=\overline{\lambda}_{n} a(x)w_{n}(x)-\overline{\lambda}_{n} w_{n}(x)u_{n}^{\rho-1}(x)\theta_{u_{n}}\leq\overline{\lambda}_{n} a(x)w_{n}(x),\eqno (3.5)$$
By the properties of operator $L$, based on (3.5), we have
$$
\|w_{n}\|_{E}\leq\| \overline{\lambda}_{n}L^{-1} (a(x)w_{n}(x))\|_{E}\leq c_{0} C ||w_{n}||_{\infty}||a||_{\infty}.
\eqno (3.6)$$
passing to the limit, then we have $\|w_{n}\|_{E}\rightarrow 0,$ which is contract with $\|w_{n}\|_{E}\equiv 1$.\hfill{$\Box$}

\noindent{\bf Proof of Theorem 1.1.}  It is obvious that any solution to (2.5) of the form $(1,u)$ yields a solution $u$ to (1.4). We
show that $\mathcal{C}$ will cross the hyperplane $1\times E$ in $\mathbb{R}^{+}\times E.$  By Lemma 2.2, Lemma 2.3 and Lemma 3.1, we just need to verify that $\lambda_{1}<1$.  In fact, it is known that the first eigenvalue $\lambda_{1}$ of the linear problem (2.1) can be given by the max-min principle using the minimum of the Rayleigh quotient, that is $$\lambda_{1}=\inf \big\{\frac{\int_{0}^{1}[u''^{2}-pu''u]dx}{\int_{0}^{1}a u^{2}dx}| \  u\in C^{4}[0,1], u\not\equiv0 \ \text{and}\  \ u(0)=u(1)=u''(0)=u''(1)=0\big\},\eqno (3.7)$$
Let us take $u$ in (3.7) as $\psi(x):=\sqrt{2}\sin (\pi x), $ then
$$\lambda_{1}\leq \frac{\int_{0}^{1}[\psi''^{2}-p\psi''\psi]dx}{\int_{0}^{1}a \psi^{2}dx}=\frac{\pi^{4}+2\pi^{2}\int_{0}^{1}p\sin^{2} (\pi x)dx}{2\int_{0}^{1}a\sin^{2} (\pi x)dx}<1.\eqno (3.8)$$\hfill{$\Box$}

\vskip 3mm

\vskip 5mm

\noindent{\bf Acknowledgements}

\noindent The authors are very grateful to the anonymous referees for their valuable
suggestions. This work was supported by the NSFC (No. 11801453).

\vskip 12mm

\centerline {\bf REFERENCES}\vskip5mm\baselineskip 0.45cm
\begin{description}
\baselineskip 15pt

\item{[1]}~  G. T. Dee, W. V. Saarloos, Bistable systems with propagating fronts leading to pattern formation, Phys. Rev. Lett. 60 (1988) 2641-2644.
\item{[2]}~  L. A. Peletier, W. C. Troy, Spatial patterns described by the extended Fisher-Kolmogorov (EFK) equation: kinks, Differential Integral Equations 8 (1995) 1279-1304.

\item{[3]}~  J. Swift, P. C. Hohenberg, Hydrodynamic fluctuations at the convective instability, Phys. Rev. A. 15 (1977) 319-328.

\item{[4]}~ R. M. Hornreich, M. Luban, S. Shtrikman, Critical behaviour at the onset of k-space instability on the ¦Ë line, Phys. Rev. Lett. 35 (1975) 1678-1681.
\item{[5]}~ M. C. Cross, P. C. Hohenberg, Pattern formation outside of equilibrium, Rev. Modern Phys. 65 (1993) 851-1112.
\item{[6]}~ Y. Chen, P. J. MacKenna, Traveling waves in a nonlinearly suspended beam: theoretical results and numerical observations, J. Differential Equations
136 (1997) 325-335.
\item{[7]}~ A. C. Lazer, P. J. MacKenna, Large-amplitude oscillations in suspension bridges: some new connections with nonlinear analysis, SIAM Rev. 32 (1990) 537-578.
\item{[8]}~ C. J. Amick, J. F. Toland, Homoclinic orbits in the dynamic phase space analogy of an elastic strut, European J. Appl. Math. 3 (1992) 97-114.
\item{[9]}~ N. N. Akhmediev, A. V. Buryak, M. Karlsson, Radiationless optical solitons with oscillating tails, Opt. Commun. 110 (1994) 540-544.
\item{[10]}~ G. J. B. van den Berg, Dynamics and Equilibria of Fourth Order Differential Equations, Selbstverl, 2000.
\item{[11]}~ A. R. Aftabizadeh; Existence and uniqueness theorem for fourth-order boundary value problems, J. Math. Anal. Appl. 116 (1986), 416-426.
\item{[12]}~ D. Bonheure, L. Sanchez, M. Tarallo and S. Terracini; Heteroclinic connections between nonconsecutive equilibria of a fourth order differential equation, Calc. Var. 17 (2003), 341-356.
\item{[13]}~ J. Chaparova, Existence and numerical aproximations of periodic solutions of semilinear fourth-order differential equations, J. Math. Anal. Appl. 273 (2002), 121-136.
\item{[14]}~ C. De Coster, C. Fabry and F. Munyamarere; Nonresonance conditions for fourth-order nonlinear boundary value problems, Internat. J. Math. Sci. 17 (1994), 725-740.
\item{[15]}~ R. Y. Ma, H. Y. Wang; On the existence of positive solutions of fourth-order ordinary differential equation, Appl. Anal. 59 (1995), 225-231.
\item{[16]}~ A. M. Micheletti, A. Pistoia; Multiplicity results for a fourth-order semilinear problem, Nonlinear Anal. 31 (1998), 895-908.
\item{[17]}~ M. A. Del Pino, R. F. Manasevich, Existence for a fourth-order nonlinear boundary problem under a two-parameter nonresonance contition, Proc. Amer. Math. Soc. 112 (1991), 81-86.
\item{[18]}~ M. R. Grossinho, L. Sanchez and S. A. Tersian, On the solvability of a boundary value problem for a fourth-order ordinary differential equation, Appl. Math. Letters 18 (2005), 439-444.
\item{[19]}~ L. A. Peletier and V. Rottschafer, Pattern selection of solutions of the Swift-Hohenberg equation, Physica D 194 (2004), 95-126.
\item{[20]}~ D. Smets, J. B. van den Berg, Homoclinic solutions for Swift-Hohenberg and suspension bridge type equations, J. Diff. Eqns. 184 (2002), 78-96.
\item{[21]}~ G. Tarantello, A note on a semilinear elliptic value problem, Differential Integral Equations, 5 (1992) 561-565.
\item{[22]}~ S. Tersian, J. Chaparova, Periodic and homoclinic solutions of extended Fisher-Kolmogorov equations, J. Math. Appl. Anal. 260 (2001), 490-506.
\item{[23]}~ G. J. B van den Berg, L. A. Peletier and W. C. Troy. Global branches of multy bump periodic
solutions of Swift-Hohenberg equation, Arch. Rational Mech. Anal. 158 (2001), 91-153.
\item{[24]}~ L. A. Peletier, W. C. Troy and V. der Vorst. Stationary solutions of a fourth-order nonlinear
diffusion equation, Differential equations, 31, 2 (1995) 301-314.
\item{[25]}~ T. Gyulov, S. Tersian, Existence of trivial and nontrivial solutions of a fourth-order ordinary differential equation, Electronic Journal of Differential Equations,  41 (2004), 1-14.
\item{[26]}~ A. L. A. de Araujo,  Periodic solutions for extended Fisher-Kolmogorov and Swift-Hohenberg equations obtained using a continuation theorem, Nonlinear Analysis 94 (2014) 100-106.
\item{[27]}~ R. Y. Ma, G. W. Dai, Periodic solutions of nonlocal semilinear fourth-order differential equations, Nonlinear Analysis 74 (2011) 5023-5029.
\item{[28]}~ P. C. Carriao, L. F. O. Faria, O. H. Miyagaki, Periodic solutions for extended Fisher-Kolmogorov and Swift-Hohenberg equations by truncature techniques, Nonlinear Anal. 67 (2007) 3076-3083.
\item{[29]}~  Z. B. Bai, H. Y. Wang, On positive solutions of some nonlinear fourth order beam equations, J. Math. Anal. Appl. 270 (2002), 357-368.
\item{[30]} A.  Cabada, J. Cid, L. Sanchez,  Positivity and lower and upper solutions for fourth order boundary value problems. Nonlinear Anal. 67(5) (2007),  1599-1612.
\item{[31]}~ J. A. Cid, D. Franco and F. Minh\'{o}s,  Positive fixed points and fourth-order equations, Bull. Lond. Math. Soc. 41 (2009), 72-78.
\item{[32]}~ P. Dr\'{a}bek, G. Holubov\'{a}, Positive and negative solutions of one-dimensional beam equation, Appl. Math. Lett. 51 (2016), 1-7.
\item{[33]}~  C. P. Gupta, Existence and uniqueness theorems for the bending of an elastic beam equation. Appl. Anal. 26(4) (1988), 289-304.
\item{[34]}~  X. L. Liu, W. T. Li, Existence and multiplicity of solutions
for fourth-order boundary value problems with parameters, J. Math.
Anal. Appl. 327 (2007), 362-375.
\item{[35]}~ Y. X. Li, A monotone iterative technique for solving the bending elastic beam equations.
Appl. Math. Comput. 217 (2010), 2200-2208.
\item{[36]}~ D. Q. Jiang, W. J. Gao, A. Y. Wan, A monotone method for constructing extremal
solutions to fourth-order periodic boundary value problems, Appl. Math. Comput. 132 (2002), 411-421.
\item{[37]}~ R. Vrabel, On the lower and upper solutions method for the problem of elastic beam with hinged ends. J. Math. Anal. Appl. 421(2) (2015), 1455-1468.
\item{[38]}~ J. R. L. Webb, G. Infante, D. Franco,
Positive solutions of nonlinear fourth-order boundary value problems
with local and non-local boundary conditions,
  Proc. Roy. Soc. Edinburgh Sect. A 138(2) (2008), 427-446.
\item{[39]}~ Z. Yang, Existence and uniqueness of positive solutions for higher order boundary value problem, Comput. Math. Appl. 54(2) (2007), 220-228.
\item{[40]}~ R. Y. Ma, Nodal solutions of boundary value problems of
fourth-order ordinary differential equations, J. Math. Anal. Appl.
319(2) (2006), 424-434.
\item{[41]}~ J. X. Wang, R. Y. Ma, S-shaped connected component for the fourth-order boundary value problem,  Boundary Value Problems,   2016) 2016:189

\item{[42]}~ P. Korman, Uniqueness and exact multiplicity of solutions for a class of fourth-order semilinear problems, Proc. Roy. Soc. Edinburgh Sect. A 134(1) (2004), 179-190.

\item{[43]}~ B. P. Rynne, Global bifurcation for $2m$th-order boundary
value problems and infinitely many solutions of superlinear
problems, J. Differential Equations 188 (2003), 461-472.

\item{[44]}~ U. Elias, Eigenvalue problems for the equation Ly+¦Ëp(x)y=0, J. Differential Equations 29 (1978) 28-57.
    
\item{[45]}~ P. H. Rabinowitz, Some global results for nonlinear eigenvalue problems, J. Funct. Anal. 7 (1971) 487-513.



\end{description}
\end{document}